\documentclass[12pt,twoside]{amsart}

\usepackage{amsopn,amssymb}
\usepackage[margin=1in]{geometry}
\usepackage[pdfauthor={David Wayne Cook II},
            pdftitle={Nauty in Macaulay2},
            pdfproducer={LaTeX with hyperref},
            pdfcreator={latex->dvips->ps2pdf},
            pdfborder={0 0 0},
            colorlinks=true]{hyperref}  

\newcommand{\st}{\; | \;}

\renewcommand{\section}[2] {\vspace{0.75\baselineskip} \noindent {\scshape #2}.}
\newtheorem*{theorem}{Theorem}
\theoremstyle{definition}
\newtheorem*{caveat}{Caveat}

\begin{document}

\title{Nauty in Macaulay2}

\author[D.\ Cook II]{David Cook II}
\address{Department of Mathematics, University of Kentucky, 715 Patterson Office Tower, Lexington, KY 40506-0027, USA}
\email{\href{mailto:dcook@ms.uky.edu}{dcook@ms.uky.edu}}
\subjclass[2010]{05C25}
\thanks{{\em Nauty} version 1.4 available at \url{http://www.ms.uky.edu/~dcook/files/Nauty.m2}; {\bf nauty} version 2.4r2.}

\begin{abstract}
    We introduce a new {\em Macaulay2} package, {\em Nauty}, which gives access to powerful methods on graphs
    provided by the software {\bf nauty} by Brendan McKay.  The primary motivation for accessing {\bf nauty} is
    to determine if two graphs are isomorphic.  We also implement methods to generate families of graphs
    restricted in various ways using tools provided with the software {\bf nauty}.
\end{abstract}

\maketitle

\section*{Introduction}
Let $G$ and $H$ be two finite, simple, undirected graphs on the common vertex set $V$ with edge sets $E(G)$ and $E(H)$, respectively.  We say
that $G$ and $H$ are {\em isomorphic} if there is a bijection $\varphi$ from $V$ to itself which preserves edges, that is, $\{u, v\} \in E(G)$
if and only if $\{\varphi(u), \varphi(v)\} \in E(H)$.  Determining whether two given graphs are isomorphic is known as the Graph Isomorphism problem.

When Garey and Johnson wrote their classic book~\cite{GJ} on the complexity of algorithms, they specified twelve problems of ambiguous complexity,
one of which was the Graph Isomorphism problem.  Unfortunately, it is still unknown if the Graph Isomorphism problem is {\bf P} or
{\bf NP}-complete.  Moreover, the problem is of such notoriety that some have even begun referring to a new complexity class, {\bf GI}, of problems
which reduce in polynomial time to the Graph Isomorphism problem~\cite{J}.  Despite this, there exists computer software which is capable of
determining whether two graphs are isomorphic in reasonable time.  One such piece of software is {\bf nauty}~\cite{N} by McKay.

The {\bf nauty} software is written in highly portable {\em C} and is designed to, above all else, compute whether two graphs are isomorphic.  
It also includes an extensive family of tools, collectively called {\bf gtools}, to generate systematic modifications of graphs, to generate
specific families of graphs, to generate random graphs, to filter a set of graphs for given properties, and to canonically relabel graphs. 
Most, if not all, of these features would be beneficial to any computer software that handles graphs.

The package {\em EdgeIdeals}~\cite{FHT} by Francisco, Hoefel, and Van Tuyl implements structures and methods for manipulating graphs
(and hypergraphs) within {\em Macaulay2}~\cite{M2}, a software system by Grayson and Stillman designed to aid in research of commutative
algebra and algebraic geometry.  We introduce a new package, {\em Nauty}, for {\em Macaulay2}, which provides an interface with
{\bf nauty}.\footnote{Throughout this note we use {\bf nauty} to refer to the software and {\em Nauty} to refer to the interface.}
Most of the aforementioned tools in {\bf gtools} are accessible through {\em Nauty}.  In particular, the methods of perhaps the greatest 
interest are {\tt areIsomorphic}, {\tt filterGraphs}, {\tt generateGraphs}, and {\tt generateRandomGraphs}.

The remainder of this note is broken in to two sections:  the first describes briefly the theoretical underpinnings of {\bf nauty} and the second
gives an example session of using {\em Nauty} along with a few useful caveats.

\section*{Canonical labellings}
In~\cite{M}, McKay describes the improved algorithms which he developed to canonically label a graph;  these algorithms are the heart
of {\bf nauty} and are summarised in~\cite{N}.  We recall briefly the theoretical ideas which make such algorithms useful.

Let $G$ be a finite, simple, undirected graph on the vertex set $V$.  We call an ordered partition $\pi = (V_1, \ldots, V_m)$
of $V$ a {\em colouring} of $G$ and call $(G, \pi)$ a {\em coloured graph}.  Given a permutation $\sigma$ of $V$, we define $\sigma(V_i)$
to be the set $\{\sigma(v) \st v \in V_i\}$, $\sigma(\pi)$ to be the colouring $(\sigma(V_1), \ldots, \sigma(V_m))$, and $\sigma(G)$
to be the graph on $V$ with edge set $E(\sigma(G))$ given by $\{\sigma(u), \sigma(v)) \in E(\sigma(G)\}$ for all $\{u, v\} \in E(G)$.
If $\sigma(\pi) = \pi$, then $\sigma$ is called {\em colour-preserving}.

If $V = [n] = \{1, \ldots, n\}$ and $\pi = (V_1, \ldots, V_m)$ is a colouring of $V$, then $c(\pi)$ is the colouring 
$(\{1, \ldots, |V_1|\}, \{|V_1| + 1, \ldots, |V_1| + |V_2|\}, \ldots, \{n - |V_m| + 1, \ldots, n\}).$   A {\em canonical labelling map}
is a function $\ell$ from the set of coloured graphs with vertex set $V$ to the set of graphs with vertex set $V$ such that, for any coloured
graph $(G, \pi)$, $\ell(G, \pi) = \tau(G)$ for some permutation $\tau$ with $\tau(\pi) = c(\pi)$ and 
$\ell(\sigma(G), \sigma(\pi)) = \ell(G, \pi)$ for every permutation $\sigma$ of $V$.

\begin{theorem}{\cite[Theorem~2.2]{M}, \cite[Theorem 1]{N}}
    Let $(G, \pi)$ and $(H, \rho)$ be coloured graphs on the common finite vertex set $V$ such that $\pi$ and $\rho$ have the same number
    of vertices in the $i^{th}$ colour class, for each $i$.  Further, let $\ell$ be a canonical labelling map for graphs with vertex set 
    $V$.  Then $\ell(G, \pi) = \ell(H, \rho)$ if and only if $\sigma(G) = H$ for some colour-preserving permutation $\sigma$.
\end{theorem}

Initially, {\bf nauty} colours a graph with a single colour and then refines this colouring using a specified vertex invariant.  Then, using 
the above theorem, {\bf nauty} can determine if two graphs are isomorphic by checking if their canonical labellings, after refinement, are the
same.  The power of {\bf nauty} is that it implements fifteen different vertex invariants (see~\cite[Section~9]{N}), each of which is more 
or less useful depending on the class of graphs being tested.

\section*{Examples}
We first load {\em Nauty}, which automatically loads {\em EdgeIdeals} for access to the {\tt Graph} class.
The {\bf nauty} software stores graphs in two different string formats, Graph6 and Sparse6 (see~\cite[Section~19]{N} for a complete
description); both formats are handled by {\em Nauty} in the method {\tt stringToGraph}.  The method {\tt graphToString} always returns
the Graph6 string-representation of the graph; e.g., the five-cycle can be represented as ``{\tt Dhc}'' and the complete graph on five
vertices is represented as ``{\tt D\~{}\{}'' in the Graph6 string-representation.

\begin{verbatim}
    i1 : needsPackage "Nauty";
    i2 : R = QQ[a..e];
    i3 : graphToString cycle R
    o3 = Dhc
    i4 : graphToString completeGraph R
    o4 = D~{
    i5 : edges stringToGraph("Dhc", R)
    o5 = {{a, b}, {b, c}, {c, d}, {a, e}, {d, e}}
\end{verbatim}

\begin{caveat}
    Due to the way graphs are constructed in {\em EdgeIdeals}, the conversion process from string format to graph format is costly.
    Thus, when many consecutive manipulations will be done through {\em Nauty}, we recommended that the graphs be left in
    string format until it is necessary to convert back to graph format.
\end{caveat}

The most powerful feature of {\em Nauty} is determining if two graphs are isomorphic.  We can demonstrate it for two given
graphs with {\tt areIsomorphic}; further we can reduce a list of graphs to be pairwise non-isomorphic with {\tt removeIsomorphs}.
We demonstrate the latter on a list of $5! = 120$ different labellings of the five-cycle, which is represented as ``{\tt Dhc}''
in the Graph6 string-representation.

\begin{verbatim}
    i6 : G = graph {{a, c}, {c, e}, {e, b}, {b, d}, {d, a}};
    i7 : areIsomorphic(cycle R, G)
    o7 = true
    i8 : removeIsomorphs apply(permutations gens R, P ->
            graphToString graph apply(5, i-> {P_i, P_((i+1)%5)}))
    o8 = {Dhc}
\end{verbatim}

{\em Nauty} also includes methods for generating all graphs on a given number of vertices, possibly with restrictions
to simple properties such as the number of edges, with the method {\tt generateGraphs}.  We generate all graphs with 
between one and nine vertices and verify the counts with the Online Encyclopedia of Integer Sequences~\cite[A000088]{OEIS}.

\begin{caveat}
    As {\bf nauty} does not handle graphs with zero vertices, {\em Nauty} will throw an error when graphs with zero
    vertices are requested or encountered.
\end{caveat}

\begin{verbatim}
    i9 : A000088 = apply(1..9, n -> #generateGraphs n)
    o9 = (1, 2, 4, 11, 34, 156, 1044, 12346, 274668)
\end{verbatim}

We can also select all graphs in a list with a given property.  To do this, we first build a filter with the method
{\tt buildGraphFilter} and then use the method {\tt filterGraphs} to select all graphs in the list which pass the filter;
see the documentation for {\tt buildGraphFilter} to see all the possible properties which can be filtered for.  First we
generate all bipartite graphs with between one and twelve vertices.  We then filter this list for forests, i.e., graphs 
without cycles.  We verify the counts with~\cite[A005195]{OEIS}.

\begin{verbatim}
    i10 : B = apply(1..12, n -> generateGraphs(n, OnlyBipartite => true));
    i11 : forestsOnly = buildGraphFilter {"NumCycles" => 0};
    i12 : A005195 = apply(B, graphs -> #filterGraphs(graphs, forestsOnly))
    o12 = (1, 2, 3, 6, 10, 20, 37, 76, 153, 329, 710, 1601)
\end{verbatim}

\begin{caveat}
    When filtering graphs for connectivity, {\em Nauty} uses the following definition of $k$-connectedness.  Specifically, 
    a $0$-connected graph is a disconnected graph and, for $k > 0$, a $k$-connected graph is a graph that can be disconnected
    by removing $k$ vertices but not by removing $k-1$ vertices.
\end{caveat}

Thus, if we wish to select only connected graphs, then we must create a filter for graphs which are {\em not} $0$-connected.
Alternatively, the method {\tt generateGraphs} has the option {\tt OnlyConnected} which forces only connected graphs
to be returned.  We demonstrate the former by filtering the list of bipartite graphs for trees, i.e., connected graphs
without cycles.  We verify the counts with~\cite[A000055]{OEIS}.

\begin{verbatim}
    i13 : treesOnly = buildGraphFilter {"NumCycles" => 0,
             "Connectivity" => 0, "NegateConnectivity" => true};
    i14 : A000055 = apply(B, graphs -> #filterGraphs(graphs, treesOnly))
    o14 = (1, 1, 1, 2, 3, 6, 11, 23, 47, 106, 235, 551)
\end{verbatim}

Last, when testing conjectures---especially when dealing with graphs on many vertices---it is nice to be able to generate
large lists of random graphs.  {\em EdgeIdeals} provides the method {\tt randomGraph}, but this method only generates
graphs with a specified number of edges and, moreover, generates only one graph at a time.  {\em Nauty} provides the
method {\tt generateRandomGraphs}, which randomly generates multiple graphs on a fixed number of vertices with a specified
edge probability.

Erd\H{o}s and R\'enyi showed in~\cite[Theorem~1]{ER} that a random graph on $n$ vertices with edge probability $\frac{(1+\epsilon)\log{n}}{n}$
is almost always connected while a graph on $n$ vertices with edge probability $\frac{(1-\epsilon)\log{n}}{n}$ is almost never connected,
at least as $n$ tends to infinity.  We demonstrate this property by generating $100$ random graphs on $n$ vertices where $\epsilon$ is
$1$ in the first case and $\frac{1}{2}$ in the second.

\begin{verbatim}
    i15 : connected = buildGraphFilter {"Connectivity" => 0, 
            "NegateConnectivity" => true};
    i16 : prob = n -> log(n)/n;
    i17 : apply(2..30, n-> #filterGraphs(
            generateRandomGraphs(n, 100, 2*(prob n)), connected))
    o17 = (70, 83, 88, 96, 94, 97, 97, 96, 97, 96, 99, 96, 99, 98, 98,
            96, 99, 96, 99, 96, 98, 97, 97, 98, 100, 99, 97, 99, 99)
    i18 : apply(2..30, n-> #filterGraphs(
            generateRandomGraphs(n, 100, (prob n)/2), connected))
    o18 = (17, 8, 3, 2, 3, 0, 1, 0, 1, 2, 3, 0, 1, 0, 0, 1, 1, 2, 0,
            2, 0, 1, 0, 0, 0, 0, 0, 0, 0)
\end{verbatim}

\section*{Acknowledgement}
We would like to thank the anonymous referee for helpful comments, especially regarding the package documentation.  We would also
like to thank the Editorial Board of the journal for many instructive comments on documentation and coding in {\em Macaulay2}.


\end{document}